\documentclass[pre,onecolumn,10pt]{revtex4}%
\usepackage{amsfonts}
\usepackage{amsmath}
\usepackage{amssymb}
\usepackage{graphicx}
\usepackage{scalefnt}
\usepackage{qtree}
\newcommand{\ignore}[1]{}
\newtheorem{theorem}{Theorem}

\newtheorem{corollary}[theorem]{Corollary}

\newtheorem{proposition}[theorem]{Proposition}

\newenvironment{proof}[1][Proof]{\noindent\textbf{#1.} }{\ \rule{0.5em}{0.5em}}
\begin{document}
\title[Diffusion on an Ising chain with kinks]{Diffusion on an Ising chain with kinks}
\author{Alioscia Hamma}
\affiliation{Massachusetts Institute of Technology, Research
Laboratory of Electronics, Cambridge MA 02139, USA}
\email{ahamma@perimeterinstitute.ca}
\author{Toufik Mansour}
\affiliation{Department of Mathematics, University of Haifa, Haifa 31905, Israel}
\email{toufik@math.haifa.ac.il}
\author{Simone Severini}
\affiliation{Institute for Quantum Computing and Department of Combinatorics \&
Optimization, University of Waterloo, Waterloo N2L 3G1, ON Canada}
\email{simoseve@gmail.com}
\begin{abstract}
We count the number of histories between the two degenerate minimum energy
configurations of the Ising model on a chain, as a function of the length $n$
and the number $d$ of kinks that appear above the critical temperature. This
is equivalent to count permutations of length $n$ avoiding certain
subsequences depending on $d$. We give explicit generating functions and
compute the asymptotics. The setting considered has a role when describing dynamics induced by quantum Hamiltonians with deconfined quasi-particles.
\end{abstract}
\maketitle

\section{Introduction}

Consider a chain with Ising variables $\sigma_{i}$ that can take the values
$\pm1$ on $n$ sites. The Ising model is described by the following action:
$S=U\sum_{i=1}^{n}(1-\sigma_{i}\sigma_{i+1})$, where $U>0$. This model has two
degenerate minimum energy configurations. Indeed, both the configurations
$\Sigma_{1,2}\equiv\{\sigma_{i}=\mp1$, for every $i\}$ minimize the energy for
this system with a value $S=0$. Every other configuration can be obtained by
flipping the $\sigma_{i}$ variables starting from either $\Sigma_{1}$ or
$\Sigma_{2}$: we see that flipping one variable costs an energy $4U$. We say that a variable is a
\emph{kink} if the adjacency variables have a different value. Flipping an
adjacent variable $\sigma_{i\pm1}$ will not change the total energy of the
system. However, if we flip a second variable $i\pm k$ with $k\geq2$,
\emph{i.e.}, we create another kink, then the cost in energy will be $8U$.
The total energy of a certain configuration $\{\sigma_{i}:i=1,...,n\}$ will be
given by $S=4U(d+1)$, where $d$ is the number of kinks, not including the
initial one.

It is obvious that the two minimum energy configurations $\Sigma_{1}$ and
$\Sigma_{2}$ are only separated by configurations that have energy $4U$. In
fact, we can start from, say $\Sigma_{1}$, pick a site $i$ and flip the
variable $\sigma_{i}$ from $-1$ to $+1$. In this way we keep staying in a
configuration with energy $4U$. Then, from there we can reach the
configuration $\Sigma_{2}$ by flipping only adjacent variables. We say that
kinks can propagate \emph{freely}, that is, without paying energy. The kink is
made of the walls of a domain with all variables of the same value. In one
dimension, making the domain larger does not change the domain walls. In
statistical mechanics, this implies that the Ising model in one dimension has
a critical temperature $T=0$. At every temperature $T>0$, some kinks will be
created, and once created, a\ kink can propagate freely and connect the two
minima (see, \emph{e.g.}, \cite{Kogut:79}). For this reason, at every $T>0$,
there is only one minimum energy configuration; there are two only when $T=0$.

We call \emph{history} the permutation on $n$ objects that indicate the
sequence of the flips connecting $\Sigma_{1}$ to $\Sigma_{2}$. So, for
instance, a history is given by the string $(i_{1},i_{2},...,i_{n})$ in which
we flip first the variable $\sigma_{1}$, then $\sigma_{2}$ and so on. As we
have pointed out, this history requires an energy $4U$. However, every time we
create a kink, we pay an energy $4U$. The energy of a history is therefore the
number $4U(d+1)$, where $d$ is number of kinks created during the entire
history, without including the first one. It is easy to compute the total
number of histories: since we must flip $n$ Ising variable, we can do it in
$n!$ possible ways, \emph{i.e.}, the number of permutations of length $n$. In
the present paper we count the total number of histories given the parameters
$n$ and $d$. We exhibit generating functions for the sequences and observe
that these numbers are asymptotically $2^{n-2d-1}(d+1)^{n}$. Our technique is
the one of generating trees (see, \emph{e.g.}, \cite{W1,W2}).

This scenario is related to all those models in condensed matter and
statistical mechanics, where two minima are separated by excitations that can
propagate freely (or quasi-freely). Even if the propagation is not free,
\emph{i.e.}, if we add a \textquotedblleft tension\textquotedblright\ to the
domains, there is a difference between propagating an excitation and creating
a new kink \cite{Kogut:79}.

Recently such a setting has become important again in the subject of
topological quantum computation \cite{Kitaev:03}. In this context, we have
quantum Hamiltonians with deconfined quasi-particles. The quasi-particles are
ends of strings. The quantum memory is encoded in the expectation value of
large, topological observables in a degenerate ground state. A crucial problem
is whether in this model quantum memory is stable for temperature higher than
absolute zero \cite{Dennis:02}. As soon as the temperature is greater than
zero, some kinks will be created with some density of probability in time, and
will start to propagate. If the total time to connect two degenerate ground
states will scale with the size of the system, the quantum memory will be
considered robust, otherwise it will be considered fragile. A precise
computation will require a knowledge about the number of all possible
histories for this process.

We will deal with a purely combinatorial framework. Let us consider a network
modeled by a chain with a finite number $n$ of nodes. In the network, the
first node is adjacent to second one; the second to the first and the third
one, and so on. To begin with, let us assume that there is a unique node of
the network having a special property. For example, the node possesses some
information, which has to be \emph{diffused} to the entire network. When the
clock ticks, the node passes the information to one its neighbors. The
process carries on, and at each time step, a single node passes the
information. It is clear that it takes exactly $n$ time steps for the
information to be transferred to the entire network.

In how many different ways this can be done? Summing over all possible initial
nodes, the answer to this question is $2^{n-1}$. This corresponds to the case
with a single kink introduced at the beginning of the process. In Figure~\ref{figS}
represented all the $2^{3-1}=4$ possible different histories, when
the number of nodes is $n=3$ (see the vertices $123$, $213$, $231$ and $321$).
The time steps of the diffusion are represented from left to right. In
this case $d=0$ because the only kink is the initial one.

Following the above description, let us now consider a variant in which we
introduce extra kinks. Again, at each time step a single node passes the
information, as we have seen before. We assume however that during the process
there are exactly $d$ nodes acquiring the information not from an adjacent
node. Additionally, these special nodes are not pairwise adjacent to each
other. The total number of ways in which the information can be propagated to
the entire network is the number of histories.

The paper is organized as follows: in Section 2 we introduce the notation and
formalize the problem; in Section 3 we state and prove the results. Section 4 is a brief conclusion.

\section{Set-up}

The graph $P_{n}=(V,E)$, where $V=\{1,2,...,n\}$ and edges
$E=\{\{1,2\},\{2,3\},...,\{n-1,n\}\}$, is called $n$\emph{-path}. We will
define below a procedure to assign labels to the vertices of $n$-paths. The
procedure will satisfy the specific conditions, giving the diffusion dynamics
which we have described in Introduction. For each integer $0\leq t\leq n$, we
define a function $f_{t}:V\longrightarrow\{+,-\}$. This gives the ordered set
$f_{t}=(f_{t}(1),f_{t}(2),...,f_{t}(n))$, specifying the image of the function
at time $t$. When $t=0$ and $t=n$, we assume that $f_{0}=(-,-,...,-)$ and
$f_{n}=(+,+,...,+)$, respectively. When $t=1$, we chose a vertex $v\in V$ and
set $f_{1}(v)=+$. In other words we designate a vertex to be a kink. At each
time step $t\geq2$, we chose a vertex $w$ such that
\begin{equation}%
\begin{tabular}
[c]{lllll}%
$f_{t}(w)=+$ & $\qquad$only if$\qquad$ & $f_{t-1}(z)=+$ & $\qquad$and$\qquad$
& $\{w,z\}\in E.$%
\end{tabular}
\ \ \ \ \label{reg}%
\end{equation}
The described rule allows us to \emph{transform }$f_{0}$ into $f_{n}$ in
exactly $n$ steps. To each one of the ordered sets $f_{0},f_{1},...,f_{n}$, we
can naturally associate a permutation $\pi$ of the set of $n$ points
$[n]=\{1,2,...,n\}$. A permutation $\pi$ is interpreted as a word of length
$n$ over the alphabet $[n]$, in which every symbol of $[n]$ appears exactly
once. Specifically, $\pi(i)=j$ if $f_{i}(j)=+$ and $f_{i-1}(j)=-$. Moreover,
by following Eq. (\ref{reg}) we have $\pi(i)=j$ if there is $k<i$, such that
$\pi(k)=j-1$ or $\pi(k)=j+1$.

Let $F_{n}$ be the set of all permutations of $[n]$ satisfying the above
conditions. The $k$-th permutation $\pi^{(k)}\in F_{n}$ is associated to the
ordered set $(f_{1}^{(k)},f_{2}^{(k)},...,f_{n}^{(k)})$. Notice that $k$ is
simply an index to distinguish between maps in $F_{n}$, being $F_{n}$ an
unordered set. An an example, the table below includes all elements of $F_{4}%
$. The left column represents $(f_{1}^{(k)},f_{2}^{(k)},...,f_{3}^{(k)})$. It
is superfluous to include $f_{0}^{(k)}$ and $f_{4}^{(k)}$, since $f_{0}%
^{(k)}=(-,-,-,-)$ and $f_{4}^{(k)}=(+,+,+,+)$ for every $k$. The right column
contains $\pi^{(k)}$. The first row is then $f_{1}^{(1)}=\left(
+,-,-,-\right)  ,f_{2}^{(1)}=\left(  +,+,-,-\right)  ,f_{3}^{(1)}=\left(
+,+,+,-\right)  $ and $\pi^{(1)}=1234$.%
\begingroup
\scalefont{0.7}
\[%
\begin{tabular}
[c]{l|l||l|l}%
$\left(  +,-,-,-\right)  ,\left(  +,+,-,-\right)  ,\left(  +,+,+,-\right)  $ &
$1234$ & $\left(  -,-,+,-\right)  ,\left(  -,+,+,-\right)  ,\left(
+,+,+,-\right)  $ & $3214$\\\hline
$\left(  -,+,-,-\right)  ,\left(  +,+,-,-\right)  ,\left(  +,+,+,-\right)  $ &
$2134$ & $\left(  -,-,+,-\right)  ,\left(  -,+,+,-\right)  ,\left(
-,+,+,+\right)  $ & $3241$\\\hline
$\left(  -,+,-,-\right)  ,\left(  -,+,+,-\right)  ,(+,+,+,-)$ & $2314$ &
$\left(  -,-,+,-\right)  ,\left(  -,-,+,+\right)  ,\left(  -,+,+,+\right)  $ &
$3421$\\\hline
$\left(  -,+,-,-\right)  ,\left(  -,+,+,-\right)  ,(-,+,+,+)$ & $2341$ &
$\left(  -,-,-,+\right)  ,\left(  -,-,+,+\right)  ,\left(  -,+,+,+\right)  $ &
$4321$%
\end{tabular}
\ \ .
\]%
\endgroup
It is simple to see that $F_{n}\varsubsetneq S_{n}$ if $n\geq3$, where $S_{n}$
denotes the set of all permutations of $[n]$. Let us consider a permutation
$\pi\notin F_{4}$. For instance, let $\pi=1324$. The intermediate steps taking
$(-,-,-,-)$ to $(+,+,+,+)$ induced by $\pi$ are $(+,-,-,-),(+,-,+,-)$ and
$(+,+,+,-)$. Notice that the second step is not in agreement with our
procedure. This is because $\pi(2)=3$ and $f_{2}(3)=+$, but $\pi(1)=1$, while,
in order to agree with the procedure, we should have $\pi(1)=2$ or $\pi(1)=4$.
The permutations in the set $S_{n}-F_{n}$ are then characterized by an extra
parameter, which we denote by $d$. This parameter specifies the number of
vertices from which the procedure can start independently. We denote by
$\pi^{(d,k)}$ the $k$-th permutation, with exactly $d$ \emph{non-adjacent
sequences} of $+$'s, inducing the functions $f_{1}^{(d,k)},f_{2}%
^{(d,k)}...,f_{n-1}^{(d,k)}$. When $d=0$, we have the above setting. In this
case, we simply wrote $\pi^{(k)}$ instead of $\pi^{(0,k)}$. The set of the
permutations of $[n]$ of the form $\pi^{(d,k)}$ is denoted by $F_{n}^{d}$.
There is a clear relation between $n$ and $d$; when $n$ is even, $d\leq n/2$;
when $n$ is odd, $d\leq(n-1)/2$. In the next section, we will find an explicit
formula for $\#F_{n}^{d}$. It is a simple fact that $S_{n}=\bigcup
\nolimits_{d}F_{n}^{d}$, every $n$. In the table below, we include the
permutations in $F_{4}^{1}\equiv S_{4}-F_{4}$. Since $\#F_{4}=8$, then
$\#F_{4}^{1}=16$.%
\begingroup
\scalefont{0.7}%
\[%
\begin{tabular}
[c]{l|l||l|l}%
$\left(  +,-,-,-\right)  ,\left(  +,+,-,-\right)  \left(  +,+,-,+\right)  $ &
$1243$ & $\left(  -,-,+,-\right)  ,\left(  +,-,+,-\right)  ,\left(
+,+,+,-\right)  $ & $3124$\\\hline
$\left(  +,-,-,-\right)  ,\left(  +,-,+,-\right)  ,\left(  +,+,+,-\right)  $ &
$1324$ & $\left(  -,-,+,-\right)  ,\left(  +,-,+,-\right)  ,\left(
+,-,+,+\right)  $ & $3142$\\\hline
$\left(  +,-,-,-\right)  ,\left(  +,-,+,-\right)  ,\left(  +,-,+,+\right)  $ &
$1342$ & $\left(  -,-,+,-\right)  ,\left(  -,-,+,+\right)  ,\left(
+,-,+,+\right)  $ & $3412$\\\hline
$\left(  +,-,-,-\right)  ,\left(  +,-,-,+\right)  ,\left(  +,+,-,+\right)  $ &
$1423$ & $\left(  -,-,-,+\right)  ,\left(  +,-,-,+\right)  ,\left(
+,+,-,+\right)  $ & $4123$\\\hline
$\left(  +,-,-,-\right)  ,\left(  +,-,-,+\right)  ,\left(  +,-,+,+\right)  $ &
$1432$ & $\left(  -,-,-,+\right)  ,\left(  +,-,-,+\right)  ,\left(
+,-,+,+\right)  $ & $4132$\\\hline
$\left(  -,+,-,-\right)  ,\left(  +,+,-,-\right)  ,\left(  +,+,-,+\right)  $ &
$2143$ & $\left(  -,-,-,+\right)  ,\left(  -,+,-,+\right)  ,\left(
+,+,-,+\right)  $ & $4213$\\\hline
$\left(  -,+,-,-\right)  ,\left(  -,+,-,+\right)  ,\left(  +,+,-,+\right)  $ &
$2413$ & $\left(  -,-,-,+\right)  ,\left(  -,+,-,+\right)  ,\left(
-,+,+,+\right)  $ & $4231$\\\hline
$\left(  -,+,-,-\right)  ,\left(  -,+,-,+\right)  ,\left(  -,+,+,+\right)  $ &
$2431$ & $\left(  -,-,-,+\right)  ,\left(  -,-,+,+\right)  ,\left(
+,-,+,+\right)  $ & $4312$%
\end{tabular}
\ \ .
\]%
\endgroup

\section{Enumeration}
\begin{figure}[htp]
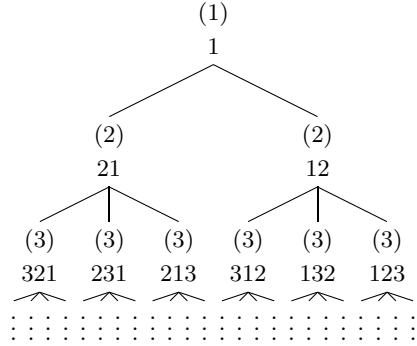

\begin{center}
{\small
\leaf{$\vdots$}\leaf{$\vdots$}\leaf{$\vdots$}\leaf{$\vdots$}\branch{4}{$\begin{array}{c}(3)\\321\end{array}$}
\leaf{$\vdots$}\leaf{$\vdots$}\leaf{$\vdots$}\leaf{$\vdots$}\branch{4}{$\begin{array}{c}(3)\\231\end{array}$}
\leaf{$\vdots$}\leaf{$\vdots$}\leaf{$\vdots$}\leaf{$\vdots$}\branch{4}{$\begin{array}{c}(3)\\213\end{array}$}
\branch{3}{$\begin{array}{c}(2)\\21\end{array}$}
\leaf{$\vdots$}\leaf{$\vdots$}\leaf{$\vdots$}\leaf{$\vdots$}\branch{4}{$\begin{array}{c}(3)\\312\end{array}$}
\leaf{$\vdots$}\leaf{$\vdots$}\leaf{$\vdots$}\leaf{$\vdots$}\branch{4}{$\begin{array}{c}(3)\\132\end{array}$}
\leaf{$\vdots$}\leaf{$\vdots$}\leaf{$\vdots$}\leaf{$\vdots$}\branch{4}{$\begin{array}{c}(3)\\123\end{array}$}
\branch{3}{$\begin{array}{c}(2)\\12\end{array}$} \branch{2}{$\begin{array}{c}(1)\\1\end{array}$}
\qobitree} \caption{Generating the permutations according to their length}\label{figS}
\end{center}
\end{figure}

A generating tree consists of a \emph{root label} and a set of
\emph{succession rules}. Following \cite{W2}, a \emph{generating
tree} is a rooted labeled tree with the property that if $v_{1}$
and $v_{2}$ are any two vertices with the same label and $l$ is
any label, then $v_{1}$ and $v_{2}$ have exactly the same number
of children with the label $l$. To specify a generating tree it
therefore suffices to specify \emph{(i)} the label of the root,
and \emph{(ii) }a set of succession rules explaining how to derive
from the label of a parent the labels of all of its children. For
example, a portion of the generating tree $\mathcal{G}$ for the
set of all unrestricted permutations is illustrated in Figure~\ref{figS}. 

The root is indexed by the permutation of length $1$. A vertex indexed by a
permutation $\pi$ of length $n$ has $n+1$ children. These are respectively
indexed by the $n+1$ permutations that can be obtained by inserting the letter
$n+1$ in the word $\pi_{1}\pi_{2}\cdots\pi_{n}$. Clearly, this tree is
isomorphic to a simpler tree, in which the root is labeled by the number $1$,
and a node labeled by $(n)$ has exactly $n+1$ children, each one
labeled by $(n+1)$:

\begin{itemize}
\item \textbf{Root}: $(1)$;

\item \textbf{Rule}: $(s)\rightsquigarrow(s+1)^{s+1}$.
\end{itemize}

\begin{figure}[h]
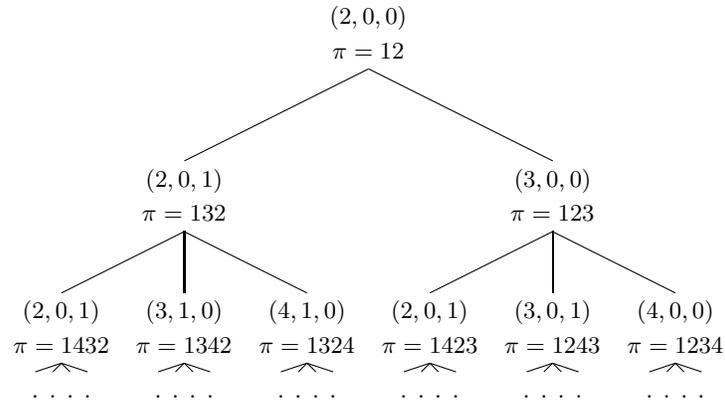

\begin{center}
{\small
\leaf{.}\leaf{.}\leaf{.}\leaf{.}\branch{4}{$\begin{array}{c}(2,0,1)\\\pi=1432\end{array}$}
\leaf{.}\leaf{.}\leaf{.}\leaf{.}\branch{4}{$\begin{array}{c}(3,1,0)\\\pi=1342\end{array}$}
\leaf{.}\leaf{.}\leaf{.}\leaf{.}\branch{4}{$\begin{array}{c}(4,1,0)\\\pi=1324\end{array}$}
\branch{3}{$\begin{array}{c}(2,0,1)\\\pi=132\end{array}$}
\leaf{.}\leaf{.}\leaf{.}\leaf{.}\branch{4}{$\begin{array}{c}(2,0,1)\\\pi=1423\end{array}$}
\leaf{.}\leaf{.}\leaf{.}\leaf{.}\branch{4}{$\begin{array}{c}(3,0,1)\\\pi=1243\end{array}$}
\leaf{.}\leaf{.}\leaf{.}\leaf{.}\branch{4}{$\begin{array}{c}(4,0,0)\\\pi=1234\end{array}$}
\branch{3}{$\begin{array}{c}(3,0,0)\\\pi=123\end{array}$}
\branch{2}{$\begin{array}{c}(2,0,0)\\\pi=12\end{array}$}
\qobitree} \caption{Generating the permutations $\pi=\pi_1\pi_2\cdots\pi_n$ with $\pi_1=1$ according to Vec$(\pi)$}\label{fig1}
\end{center}
\end{figure}

In order to enumerate the number of permutations in $F_{n}^{d}$, we define the
generating tree $\mathcal{G}$, in which we replace each permutation $\pi
=\pi_{1}\pi_{2}\cdots\pi_{n}$ with $\pi_1=1$ by the vector Vec$(\pi)=(j,k,r)$, where $\pi_{j}=n$, $\pi\in F_{n}^{k}$ and
\[
r=\left\{
\begin{array}
[c]{ll}%
1, & \pi_{n}^{-1}<\pi_{n-1}^{-1};\\
0, & \pi_{n}^{-1}>\pi_{n-1}^{-1},
\end{array}
\right.
\]
as it is shown in Figure~\ref{fig1}.

\begin{proposition}
\label{pro1} The generating tree $\mathcal{G}$ is given by

\begin{itemize}
\item \textbf{Roots: }$(2,0,0)$ and $(1,0,1)$;

\item \textbf{Rules: }

$(j,k,0)\rightsquigarrow(1,k+1,1)\cdots(j,k+1,1)(j+1,k,0)\cdots(n+1,k,0)$;

$(j,k,1)\rightsquigarrow(1,k,1)\cdots(j,k,1)(j+1,k,0)\cdots(n+1,k,0)$.
\end{itemize}
\end{proposition}

\begin{proof}
From the definitions, the permutations $12$ and $21$ are associated to the
vectors $(2,0,0)$ and $(1,0,1)$, respectively. Let $\pi=1\pi_{2}\cdots\pi_{n}$
be a permutation in $S_{n}$ associated to the vector $(j,k,0)$. By inserting
the letters $n+1$ between position $i-1$ and position $i$, we obtain a new
permutation $\pi=1\pi_{2}\cdots\pi_{i-1}(n+1)\pi_{i}\cdots\pi_{n}\in S_{n+1}$
associated to the vector $(i,k+1,1)$, for all $i=2,3,\ldots,j$, and the
vector $(i,k,0)$, for all $i=j+1,j+2,\ldots,n+1$. Let $\pi=1\pi_{2}\cdots
\pi_{n}$ be a permutation in $S_{n}$ associated to the vector $(j,k,1)$. By
inserting the letters $n+1$ between position $i-1$ and position $i$, we obtain
a new permutation $\pi=1\pi_{2}\cdots\pi_{i-1}(n+1)\pi_{i}\cdots\pi_{n}\in
S_{n+1}$ associated to the vector $(i,k,1)$, for all $i=1,2,\ldots,j$, and the
vector $\ (i,k,0)$, for all $i=j+1,j+2,\ldots,n+1$. This completes the proof.
\end{proof}

Let $G_{0}(t;u,v)$ and $G_{1}(t;u,v)$ be the associated generating functions:%
\begin{align*}
G_{r}(t;u,v)  &  =\sum_{n\geq2}t^{n}\sum_{j\geq2}\sum_{k\geq0}\sum_{\pi\in
S_{n}\,:\text{Vec}(\pi)=(j,k,r)}u^{j}v^{k}\\
&  =\sum_{n\geq2}t^{n}\sum_{j\geq2}\sum_{k\geq0}g_{n,j,k,r}u^{j}v^{k},
\end{align*}
where $r=0,1$. Alternatively, this series counts the vertices of the tree
$\mathcal{G}$ by their height and their label (the root being at height $2$).

\begin{proposition}
\label{pro2}The generating function for the number of permutations
$\pi=\pi_{1}\cdots\pi_{n}$ in $S_{n}$, $n\geq2$, according to the statistics
$j,k$, where $\pi_{j}=n$ and $\pi\in F_{n}^{k}$, is given by
$G_{0}(t;u,v)+G_{1}(t;u,v)$ such that the generating functions $G_{0}(t;u,v)$ and
$G_{1}(t;u,v)$ satisfy%
\[
G_{0}(t;u,v)=t^{2}u^{2}+\frac{tu}{1-u}(G_{0}(t;u,v)-uG_{0}(tu;1,v))+\frac
{tu}{1-u}(G_{1}(t;u,v)-uG_{1}(tu;1,v))
\]
and%
\[
G_{1}(t;u,v)=t^{2}u+\frac{tuv}{1-u}(G_{0}(t;1,v)-G_{0}(t;u,v))+\frac{tu}%
{1-u}(G_{1}(t;1,v)-G_{1}(t;u,v)).
\]
\end{proposition}
\begin{proof}
Since $g_{n,j,k,r}$ is the number of vertices in the generating tree $\mathcal{G}$, as it was shown in
Proposition~\ref{pro1}, with label $(j,k,r)$, we obtain that 
\begin{align*}
g_{n+1,j,k,0}&=\sum_{i=1}^{j-1} (g_{n,i,k,0}+g_{n,i,k,1}),\\
g_{n+1,j,k,1}&=\sum_{i=j}^n (g_{n,i,k-1,0}+g_{n,i,k,1}).
\end{align*}
Let $g_{n+1,r}(u,v)=\sum_{j\geq1}\sum_{k\geq0}g_{n+1,j,k,r}u^jv^k$ for $r=0,1$. Rewriting the above recurrence relations in terms of $g_{n,r}(u,v)$, we obtain that
\begin{align*}
g_{n+1,0}(u,v)&=\sum_{j\geq1}\sum_{k\geq1}u^jv^k\left(\sum_{i=1}^{j-1}(g_{n,i,k,0}+g_{n,i,k,1})\right),\\
g_{n+1,1}(u,v)&=\sum_{j\geq1}\sum_{k\geq1}u^jv^k\left(\sum_{i=j}^n (g_{n,i,k-1,0}+g_{n,i,k,1})\right),
\end{align*}
which is equivalent to 
\begin{align*}
g_{n+1,0}(u,v)&=\sum_{j\geq1}\sum_{k\geq0}g_{n,j,k,0}v^{k}(u^{j+1}+\cdots+u^{n+1})+\sum_{j\geq1}\sum_{k\geq0}g_{n,j,k,1}v^k(u^{j+1}+\cdots+u^{n+1}),\\
g_{n+1,1}(u,v)&=\sum_{j\geq1}\sum_{k\geq0}g_{n,j,k,0}v^{k+1}(u^1+\cdots+u^j)+\sum_{j\geq1}\sum_{k\geq0}g_{n,j,k,1}v^k(u^1+\cdots+u^j).
\end{align*}
Multiplying by $t^{n+1}$ and summing over $n\geq2$ together with using the initial conditions $g_{2,0}(u,v)=u^2$ and $g_{2,1}(u,v)=u$, we obtain that
\begin{align*}
G_{0}(t;u,v)  &  =t^{2}u^{2}+t\sum\limits_{n\geq2}\sum\limits_{j\geq1}%
\sum\limits_{k\geq0}g_{n,j,k,0}t^{n}v^{k}(u^{j+1}+\cdots+u^{n+1}%
)+t\sum\limits_{n\geq2}\sum\limits_{j\geq1}\sum\limits_{k\geq0}g_{n,j,k,1}%
t^{n}v^{k}(u^{j+1}+\cdots+u^{n+1})\\
&  =t^{2}u^{2}+\frac{tu}{1-u}\sum\limits_{n\geq2}\sum\limits_{j\geq1}%
\sum\limits_{k\geq0}g_{n,j,k,0}t^{n}v^{k}(u^{j}-u^{n+1})+\frac{tu}{1-u}%
\sum\limits_{n\geq2}\sum\limits_{j\geq1}\sum\limits_{k\geq0}g_{n,j,k,1}%
t^{n}v^{k}(u^{j}-u^{n+1}),
\end{align*}
and%
\begin{align*}
G_{1}(t;u,v)  &  =t^{2}u+t\sum\limits_{n\geq2}\sum\limits_{j\geq1}%
\sum\limits_{k\geq0}g_{n,j,k,0}t^{n}v^{k+1}(u^{1}+\cdots+u^{j})+t\sum
\limits_{n\geq2}\sum\limits_{j\geq1}\sum\limits_{k\geq0}g_{n,j,k,1}t^{n}%
v^{k}(u^{1}+\cdots+u^{j})\\
&  =t^{2}u+\frac{tuv}{1-u}\sum\limits_{n\geq2}\sum\limits_{j\geq1}%
\sum\limits_{k\geq0}g_{n,j,k,0}t^{n}v^{k}(1-u^{j})+\frac{tu}{1-u}%
\sum\limits_{n\geq2}\sum\limits_{j\geq1}\sum\limits_{k\geq0}g_{n,j,k,1}%
t^{n}v^{k}(1-u^{j}).
\end{align*}
From the definitions of the generating functions $G_0(t;u,v)$ and
$G_1(t;u,v)$, we get that these functions satisfy the
following functional equations:
\begin{align*}
G_{0}(t;u,v)  &  =t^{2}u^{2}+\frac{tu}{1-u}(G_{0}(t;u,v)-uG_{0}(tu;1,v))+\frac{tu}%
{1-u}(G_{1}(t;u,v)-uG_{1}(tu;1,v)),
\end{align*}
and
\begin{align*}
G_{1}(t;u,v)  &  =t^{2}u+\frac{tuv}{1-u}(G_{0}(t;1,v)-G_{0}(t;u,v))+\frac{tu}{1-u}%
(G_{1}(t;1,v)-G_{1}(t;u,v)),
\end{align*}
as claimed.
\end{proof}

\begin{theorem}
\label{th01}The generating function
\[
G(t,v)=\sum_{n\geq2}\sum_{d\geq0}\#F_{n}^{d}v^{d}t^{n}%
\]
is given by%
\[
G(t,v)=\sum_{j\geq0}\frac{4t^{2}\sqrt{1-v}(1-(1+2j)\sqrt{1-v}t-t(1-v))v^{j}%
}{(1+\sqrt{1-v})^{1+2j}(1-2jt\sqrt{1-v})^{2}(1-2(j+1)t\sqrt{1-v})^{2}}.
\]
\end{theorem}
\begin{proof} Denote the generating functions
$G_{0}(t;u,v)+G_{1}(t;u,v)$ and $vG_{0}(t;u,v)+G_{1}(t;u,v)$ by
$G(t;u,v)$ and $H(t;u,v)$, respectively. Then,
Proposition~\ref{pro2} gives that%
\[\begin{array}{l}
\left(  1-\frac{t^{2}u^{2}(1-v)}{(1-u)^{2}}\right)  G(t;u,v)\\
\qquad\qquad\qquad+\frac{u^{2}%
t}{(1-u)^{2}}(1-u+ut(1-v))G(ut;1,v)-\frac{ut}{1-u}H(t;1,v)+\frac{ut^{2}%
(u^{2}-1-u^{2}t(1-v))}{1-u}=0
\end{array}
\]
and%
\[
\begin{array}{l}
\left(  1-\frac{t^{2}u^{2}(1-v)}{(1-u)^{2}}\right)  H(t;u,v)\\
\qquad\qquad\qquad+\frac{u^{2}%
vt}{1-u}G(ut;1,v)-\frac{ut}{(1-u)^{2}}(1-u-ut(1-v))H(t;1,v)+\frac{ut^{2}}%
{1-u}(u-1-vu+vu^{2}+ut(1-v))=0.\end{array}
\]
This type of functional equations can be solved systematically using
the kernel method (see \cite{HM} and references therein). In this
case, if we substitute $u=p=\frac{1}{1+t\sqrt{1-v}}$ and
$u=q=\frac{1}{1-t\sqrt{1-v}}$, we obtain that%
\[
H(t;1,v)=\frac{1+\sqrt{1-v}}{1+\sqrt{1-v}t}G(\frac{t}{1+\sqrt{1-v}%
t};1,v)+\frac{t^{2}\sqrt{1-v}(2+\sqrt{1-v}+t\sqrt{1-v})}{(1+\sqrt{1-v}t)^{2}}%
\]
and%
\[
H(t;1,v)=\frac{1-\sqrt{1-v}}{1-\sqrt{1-v}t}G(\frac{t}{1-\sqrt{1-v}%
t};1,v)+\frac{t^{2}\sqrt{1-v}(2-\sqrt{1-v}-t\sqrt{1-v})}{(1-\sqrt{1-v}t)^{2}%
}.
\]
Hence, the generating function $G(t,v)=G(t;1,v)$ satisfies%
\begin{align*}
&  \frac{1+\sqrt{1-v}}{1+\sqrt{1-v}t}G(\frac{t}{1+\sqrt{1-v}t},v)+\frac
{t^{2}\sqrt{1-v}(2+\sqrt{1-v}+t\sqrt{1-v})}{(1+\sqrt{1-v}t)^{2}}\\
&  =\frac{1-\sqrt{1-v}}{1-\sqrt{1-v}t}G(\frac{t}{1-\sqrt{1-v}t},v)+\frac
{t^{2}\sqrt{1-v}(2-\sqrt{1-v}-t\sqrt{1-v})}{(1-\sqrt{1-v}t)^{2}}%
\end{align*}
which is equivalent to%
\[
G(t,v)=\frac{4t^{2}\sqrt{1-v}}{(1-2t\sqrt{1-v})^{2}(1+\sqrt{1-v})}%
(1-t\sqrt{1-v}-t(1-v))+\frac{v}{(1+\sqrt{1-v})^{2}(1-2t\sqrt{1-v})}G\left(
\frac{t}{1-2t\sqrt{1-v}},v\right)  .
\]
Applying the above functional equation an infinite number of times,
we can find an explicit formula for the generating function
$G(t,v)$,as requested.
\end{proof}

Let $h_{n}(v)$ be the polynomial $\sum_{d=0}^{n}\#F_{n}^{d}v^{d}$.
For $n=2,3,\ldots,10$, Theorem \ref{th01} gives
\begin{align*}
h_{2}(v) &  =2,\\
h_{3}(v)&=4+2v,\\
h_{4}(v)&=8+16v,\\
h_{5}(v) &  =16+88v+16v^{2},\\
h_{6}(v)&=32+416v+272v^{2},\\
h_{7}(v) &  =64+1824v+2880v^{2}+272v^{3},\\
h_{8}(v)&=128+7680v+24576v^{2}+7936v^{3},\\
h_{9}(v) &  =256+31616v+185856v^{2}+137216v^{3}+7936v^{4},\\
h_{10}(v)&=512+128512v+1304832v^{2}+1841152v^{3}+353792v^{4}.
\end{align*}
On the other hand, Theorem \ref{th01} can be used to obtain explicit
formula for the generating functions
$h^{d}(t)=\sum_{n\geq2}\#F_{n}^{d}t^{n}$. For example, Theorem
\ref{th01}, with $v=0$, gives
$h^{0}(t)=\frac{2t^{2}}{1-2t}=\sum_{n\geq 2}2^{n-1}t^{n}$, which
implies that $\#F_{n}^{0}=2^{n-1}$, for all $n\geq2$. More
generally, in order to obtain an explicit formula for the generating
function $h^{d}(t)$, on the basis of the theorem, it is enough to
take the derivative of the sum
\[
\sum_{j=0}^{d}\frac{2t^{2}v^{j}\sqrt{1-v}}{(1-2jt\sqrt{1-v})(1-2(j+1)t\sqrt
{1-v})(1+\sqrt{1-v})^{2j+1}}%
\]
exactly $d$ times and then substitute $v=0$. Applying Theorem~\ref{th01} for $d=0$, that is looking for the free coefficient of $v$ in the generating function $G(t,v)$ (in this case also $j=0$), we get that
\[h^0(t)=\frac{2t^2}{1-2t}.\]
By considering the coefficient of $v^j$, $j=1,2,3$, in the generating function $G(t,v)$ (see Theorem~\ref{th01}), we obtain an explicit formula for the generating function $h^j(t)$, and we can state the following result.

\begin{corollary}
We have
\begin{align*}
h^{0}(t)  &  =\frac{2t^{2}}{1-2t},\qquad h^{1}(t)=\frac{2t^{3}}{(1-2t)^{2}%
(1-4t)},\\
h^{2}(t)  &  =\frac{16t^{5}(1-3t)}{(1-2t)^{3}(1-4t)^{2}(1-6t)},\qquad
h^{3}(t)=\frac{16t^{7}(17-184t+636t^{2}-720t^{3})}{(1-2t)^{4}(1-4t)^{3}%
(1-6t)^{2}(1-8t)}.
\end{align*}
\end{corollary}

This gives an exact formula for $\#F_{n}^{d}$, when $d=0,1,2,3$ and
$n\geq1,3,5,7$, respectively%
\begin{align*}
\#F_{n}^{0} &  =2^{n-1};\\
\#F_{n}^{1}&=2^{n-2}%
(2^{n-1}-n);\\
\#F_{n}^{2} &  =\frac{1}{32}6^{n}-\frac{1}{16}4^{n}\left(  n-1\right)
+\frac{1}{32}2^{n}\left(  2n^{2}-4n-1\right)  \\
\#F_{n}^{3} &  =\frac{1}{128}8^{n}-\frac{1}{64}6^{n}\left(
n-2\right) +\frac{1}{64}4^{n}\left(  n^{2}-4n+2\right)
-\frac{1}{192}2^{n}\left( 12n^{2}-13n+2n^{3}+6\right).
\end{align*}

We now focus on the asymptotic behavior of $h_{d}$. Two sequences $a_{n}$ and
$b_{n}$ are said to be \emph{asymptotically equivalent} as $n\rightarrow
\infty$ if $\lim_{n\rightarrow\infty}\frac{a_{n}}{b_{n}}=1$. We denote this
fact by $a_{n}\sim b_{n}$.

\begin{theorem}
For any fixed $d\geq0$,
\[
\#F_{n}^{d}\sim2^{n-2d-1}(d+1)^{n}.
\]
\end{theorem}
\begin{proof}
Theorem \ref{th01} shows that the smallest positive pole of the generating
function $h^{d}(t)$ is $t^{\ast}=\frac{1}{2(d+1)}$ of order $1$. Thus, we have
$\#F_{n}^{d}\sim c_{d}2^{n}(d+1)^{n}$. Concerning the constant $c_{d}$, we can
write%
\begin{align*}
c_{d}  &  =\lim\limits_{t\rightarrow t^{\ast}}(1-t/t^{\ast})\lim
\limits_{v\rightarrow0}\frac{1}{d!}\frac{\partial^{d}}{\partial v^{d}}\left(
\sum_{j\geq0}\frac{4t^{2}\sqrt{1-v}(1-(1+2j)\sqrt{1-v}t-t(1-v))v^{j}}%
{(1+\sqrt{1-v})^{1+2j}(1-2jt\sqrt{1-v})^{2}(1-2(j+1)t\sqrt{1-v})^{2}}\right)
\\
&  =\lim\limits_{t\rightarrow t^{\ast}}(1-t/t^{\ast})\lim\limits_{v\rightarrow
0}\left(  \sum_{j=0}^{d}\frac{\partial^{d-j}}{\partial v^{d-j}}\frac
{4j!{\binom{{d}}{{j}}}t^{2}\sqrt{1-v}(1-(1+2j)\sqrt{1-v}t-t(1-v))}%
{d!(1+\sqrt{1-v})^{1+2j}(1-2jt\sqrt{1-v})^{2}(1-2(j+1)t\sqrt{1-v})^{2}}\right)
\\
&  =\lim\limits_{t\rightarrow t^{\ast}}\lim\limits_{v\rightarrow0}\frac
{4t^{2}(1-t/t^{\ast})\sqrt{1-v}(1-(1+2d)\sqrt{1-v}t-t(1-v))}{(1+\sqrt
{1-v})^{1+2d}(1-2dt\sqrt{1-v})^{2}(1-2(d+1)t\sqrt{1-v})^{2}}\\
&  =\lim\limits_{t\rightarrow t^{\ast}}\frac{4t^{2}(1-t/t^{\ast})^{2}%
}{2^{1+2d}(1-2dt)^{2}(1-t/t^{\ast})^{2}}\\
&  =\frac{4(t^{\ast})^{2}}{2^{1+2d}(1-2d/t^{\ast})^{2}}\\
&  =\frac{1}{2^{2d+1}},
\end{align*}
as claimed.
\end{proof}

\section{Conclusions}
In this paper, we have counted the number $\#F_{n}^{d}$ of
histories that allow to connect a ground state of an Ising chain
with $n$ variables to the other one, having specified the number
$d$ of domain walls (or kinks) created in every history. We have
given the generating functions and studied the asymptotic
behavior.

The solution of this problem is important in statistical
mechanics, to solve the random walk with kinks, and in quantum
computing, to investigate the robustness of topological quantum
memory at finite temperature. There is a number of open problems.
A first problem would be to count the number of histories to flip
the Ising variables along a non contractible loop in a torus. A
second open problem would be to introduce the possibility that at
a given click no variable is flipped or that a variables is
flipped back. Eventually, the goal would be is to assign
probabilities to the three different events: flipping an adjacent
variable, creating a kink, no flip, and computing the probability
of making a loop around the torus in a given number of steps.


\end{document}